\documentclass[10pt,leqno]{amsart} %%%%leqno es para numeraci\'{o}n de
   \usepackage[centertags]{amsmath}
   \usepackage{amsfonts}
   \usepackage{amsmath}
   \usepackage{amssymb}
   \usepackage{amsthm}
   \usepackage{newlfont}
   \usepackage[latin1]{inputenc}%%%%para que reconozca las tildes.
\usepackage{xcolor}

\allowdisplaybreaks[3]
\theoremstyle{plain}
\newtheorem{Theo}{Theorem}

\newtheorem{Cor}[Theo]{Corollary}
\newtheorem{Lem}[Theo]{Lemma}
\newtheorem{TheoA}{Theorem}

\theoremstyle{definition}

\theoremstyle{remark}

\newcommand{\ZZ}{\mathbb{Z}}
\newcommand{\NN}{\mathbb{N}}

\numberwithin{equation}{section}

\newcommand\G{{\mathcal G}}
\newcommand\Hh{{\mathcal H}}

\newcommand\Be{\mathfrak{b}}

\newcommand\sign{\operatorname{sign}}

\newcommand*\pFq{\begingroup
        \dopFq
}
\catcode`\,12
\def\dopFq#1#2#3#4#5{%
        {}_{#1}F_{#2}\biggl(\genfrac..{0pt}{}{#3}{#4};#5\biggr)%
        \endgroup
}

\title{Generalized Bell polynomials}
\author{Antonio J. Dur\'an}
\address{Departamento de An\'a\-li\-sis Mate\-m\'a\-ti\-co and IMUS,
        Universidad de Sevilla,
        41080 Sevilla, Spain}
\email{duran@us.es}
   \date{}

   \thanks{This research was partially supported by PID2021-124332NB-C21
(Minis\-te\-rio de Cien\-cia e Inno\-va\-ci\'on and Feder Funds (European Union)), and
FQM-262 (Jun\-ta de Anda\-lu\-c\'ia).}


\begin{document}
   \maketitle

\begin{abstract}
In this paper, generalized Bell polynomials $(\Be_n^\phi)_n$ associated to a sequence of real numbers $\phi=(\phi_i)_{i=1}^\infty$ are introduced. Bell polynomials correspond to $\phi_i=0$, $i\ge 1$. We prove that when $\phi_i\ge 0$, $i\ge 1$: (a) the zeros of the generalized Bell polynomial $\Be_n^\phi$ are simple, real and non positive; (b) the zeros of $\Be_{n+1}^\phi$ interlace the zeros of $\Be_n^\phi$; (c) the zeros are decreasing functions of the parameters $\phi_i$. We find a hypergeometric representation for the generalized Bell polynomials. As a consequence, it is proved that the class of all generalized Bell polynomials is actually the same class as that of  all Laguerre multiple polynomials of the first kind.
\end{abstract}

\section{Introduction and results}
The Bell polynomials $(\Be_n)_n$ are defined by
\begin{equation}\label{ben}
\Be_n(x)=\sum_{j=0}^nS(n,j)x^j,\quad n\ge 0,
\end{equation}
where $S(n,j)$ are the Stirling numbers of the second kind
\begin{equation}\label{stir}
S(n,j)=\sum_{i=0}^j\frac{(-1)^{j-i}i^n}{(j-i)!\, i!},\quad 0\le j\le n,
\end{equation}
so that
\begin{equation}\label{comp0}
x^n=\sum_{j=0}^nS(n,j)x(x-1)\cdots (x-j+1)
\end{equation}
(see \cite{boy} for a review of Bell polynomials).

The Bell polynomials can be recursively defined as follows: for $n\ge 0$,
\begin{equation}\label{rbel}
\Be_{n+1}(x)=x\left(1+\frac{d}{dx}\right)\Be_n(x),
\end{equation}
and $\Be_0=1$.

Bell polynomials are sometimes called as Touchard or exponential polynomials and were studied by Ramanujan (see \cite[Chapter 3]{ber}) in his notebooks before they were introduced by Touchard \cite{tou} and Bell \cite{bel}. They are also called single-variable Bell polynomials to distinguish them from the $k$-variable Bell polynomials $\Be_n(x_1,\dots, x_k)$ (the single-variable case being a particular instance of the $k$-variable one:  $\Be_n(x)=\Be_n(x,\dots, x)$).

It was proved by Harper \cite{har} that the Bell polynomials have simple and non positive zeros, and that the negative zeros of $\Be_{n+1}$ interlace the zeros of $\Be_n$.

The interlacing property is defined as follows: A $k$-set is a finite set of $k$ real numbers arranged in increasing or decreasing order.
We then say that the $k$-set $U$ interlaces the $\kappa$-set $V$ if
$$
\begin{cases} \mbox{between two consecutive elements of $U$}&\\ \quad\quad \mbox{there exists one and only one element of $V,$ and}&\\\mbox{either $k=\kappa+1$, or
$k=\kappa$ and $\max U<\max V.$}\end{cases}
$$
Let us note that, according to our definition, when $k=\kappa$, the interlacing property is not symmetric.

The Bell polynomials have been mainly considered in combinatorics although have applications in probability and some other areas. For instance, the Bell polynomials can also be defined as
the moments of the Poisson distribution with expected value $x$:
$$
\Be_n(x)=\int_0^\infty t^ndw_x(t),
$$
where $w_x$ is defined by
$$
w_x=e^{-x}\sum_{j=0}^\infty \frac{x^j}{j!}\delta_j
$$
(let us note that $w_x$ can also be seen as the orthogonalizing weight for the Charlier polynomials).

The purpose of this paper is to construct generalized Bell polynomials associated to a given sequence of real numbers. Given a sequence $\phi=(\phi_i)_{i\ge 1}$ of real numbers, we define the generalized Bell polynomials $(\Be_n^\phi)_n$ as follows. For $n\ge0$, set
\begin{equation}\label{dpm}
\Phi_i^n\equiv \Phi_i^n(\phi)=\begin{cases} 1, &i=0,\\
\displaystyle \sum_{1\le j_1<\cdots<j_i\le n}\phi_{j_1}\cdots\phi_{j_i}, &1\le i\le n,\end{cases}
\end{equation}
so that
$$
\prod_{i=1}^n(x+\phi_i)=\sum_{j=0}^n\Phi_{n-j}^nx^j.
$$

The $n$-th generalized Bell polynomial $\Be_n^\phi$ is defined by
\begin{equation}\label{dpbg}
\Be_n^\phi(x)=\sum_{j=0}^n\Phi_{n-j}^n\Be_j(x).
\end{equation}
Notice that the polynomial $\Be_n^\phi$ has degree $n$, is monic, and only depends on the numbers $\phi_1,\dots, \phi_n$. Moreover, $\Be_n^\phi$ is a symmetric function of $\phi_1,\dots, \phi_n$.

The Bell polynomials correspond to the case $\phi_i=0$, $i\ge 1$. The case $\phi_i=r\not =0$, $i\ge 1$, are the $r$-Bell polynomials, which have also been studied, mainly when $r$ is a positive integer, because of their combinatorial interest (see, for instance, \cite{mez,MeCo,MeRa}). When $\phi_i=0$, for $i\ge K$, the generalized Bell polynomials reduce
to a linear combination of $K$ consecutive Bell polynomials.

The content of the paper is as follows.

In Section \ref{pre}, we study some basic properties of the generalized Bell polynomials. In particular, we establish
the following recurrence relation: for $n\ge 0$,
\begin{equation}\label{rrgp}
\Be_{n+1}^\phi(x)=x\left(1+\frac{d}{dx}\right)\Be_n^{\phi}(x)+\phi_{n+1}\Be_n^{\phi}(x),
\end{equation}
and $\Be_0^\phi=1$.

We also prove that if the numbers $\rho_{n,j}^\phi$, $0\le j\le n$, are defined from the equation
\begin{equation}\label{comp1}
\prod_{i=1}^n(x+\phi_i)=\sum_{j=0}^n\rho_{n,j}^\phi x(x-1)\cdots (x-j+1),
\end{equation}
then
\begin{equation}\label{cff}
\Be_n^\phi(x)=\sum_{j=0}^n\rho_{n,j}^\phi x^j
\end{equation}
(compare (\ref{ben}), (\ref{comp0}) with (\ref{cff}), (\ref{comp1})).
As a consequence, the generalized Bell polynomials can also be defined as the generalized moments of the Poisson distribution with expected value $x$:
$$
\Be_n^\phi(x)=\int_0^\infty (t+\phi_1)\cdots (t+\phi_n)dw_x(t),\quad n\ge 0.
$$

Using (\ref{cff}), we find the following hypergeometric representation for the generalized Bell polynomials:

\begin{equation}\label{hypg}
e^{x}\Be_n^\phi(x)=\sum_{j=0}^\infty \frac{\prod_{i=1}^n(j+\phi_i)}{j!}x^j= \left(\prod_{i=1}^n\phi_i\right)\pFq{n}{n}{1+\phi_1,\dots,1+\phi_n}{\phi_1,\dots,\phi_n}{x},
\end{equation}
where for the second identity we assume $\phi_i\not =0,-1,-2,\dots $.

The hypergeometric representation (\ref{hypg}) also established a (somehow unexpected) connection between the generalized Bell polynomials and the so-called Laguerre multiple polynomials of the first kind (see \cite{Ism}, p. 627, or also \cite{aba,newa}).
Given a multi-index $\vec{n}=(n_1,\dots, n_q)\in \NN^q$ and parameters $\vec{\alpha}=(\alpha_1,\dots, \alpha_q)$, $\alpha_i-\alpha_j\not \in \ZZ$, the Laguerre multiple polynomials of the first kind are defined by the Rodrigues formula
$$
L_{\vec{n}}^{\vec{\alpha}}(x)=(-1)^{|\vec{n}|}e^{x}\prod_{j=1}^q\left(x^{-\alpha_j}\frac{d^{n_j}}{dx^{n_j}}x^{n_j+\alpha_j}\right)e^{-x},
$$
where $|\vec{n}|=\sum_{j=1}^qn_j$.
They satisfy the orthogonality conditions ($\alpha_j>-1$)
$$
\int_0^\infty L_{\vec{n}}^{\vec{\alpha}}(x)x^{\alpha_j}e^{-x}x^k=0,\quad k=0,1,\dots, n_j-1,
$$
for $j=1,2,\dots, q$ (so, according to the usual definition of multiple orthogonal polynomials, see \cite{MW}, they are the \textit{type II} multiple Laguerre polynomials of the first kind).
The hypergeometric representation (\ref{hypg}) implies that the class of all generalized Bell polynomials is the same class as that of all Laguerre multiple polynomials of the first kind. More precisely,
\begin{equation}\label{cgblm}
L_{\vec{n}}^{\vec{\alpha}}(-x)=(-1)^{|\vec{n}|}\Be_{|\vec{n}|}^{\phi^{{\vec \alpha},\vec{n}}}(x),
\end{equation}
where
\begin{equation}\label{anph}
\phi^{\vec \alpha,\vec{n}}_i=\begin{cases}\alpha_1+i,&i=1,\dots,n_1,\\
\alpha_2+i-n_1,&i=n_1+1,\dots,n_1+n_2,\\
\cdots &\\
\alpha_q+i-(n_1+\cdots+n_{q-1}),&i=n_1+\cdots+n_{q-1}+1,\dots,|\vec{n}|.
\end{cases}
\end{equation}
Let us notice that the sequence $\phi^{{\vec \alpha},\vec{n}}$ depends on both the parameters $\vec{\alpha}$ and the multi-index $\vec{n}$.

In particular, for the Laguerre polynomials, we have
$$
L_n^{\alpha}(-x)=\frac{1}{n!}\Be_n^{\phi^\alpha}(x),
$$
where now the sequence $\phi^\alpha_i=\alpha+i$ only depends on $\alpha$.

In Section \ref{zer}, as one of the main result of this paper, we prove that when $\phi_i\ge 0$, the zeros of the generalized Bell polynomials behave very nicely.

\begin{Theo}\label{pcgbi} Assume that $\phi_i\ge 0$.
\begin{enumerate}
\item The polynomial $\Be_n^\phi$ has $n$ simple and non positive zeros, and the negative zeros of $\Be_{n+1}^\phi$ interlace the negative zeros of $\Be_{n}^{\phi}$.
\item We denote $\zeta_k(n,\phi)$, $1\le k\le n$, for the $k$-th zero of $\Be_n^\phi$, arranging the zeros in increasing order, so that $\zeta_1(n,\phi)$ is the leftmost zero (to simplify the notation and when the context allows it we sometimes will write $\zeta_k$, $\zeta_k(n)$ or $\zeta_k(\phi)$). Then $\zeta_k(\phi)$ is a decreasing function of $\phi$ (where we say that $\phi\prec \psi$ if $\phi_i\le \psi_i$, $i\ge 1$, and there exists $i_0$ such that $\phi_{i_0}<\psi_{i_0}$).
\end{enumerate}
\end{Theo}

As a consequence of the connection with the Laguerre multiple polynomials of the first kind, we deduce that the zeros of these polynomials are increasing functions of the parameters $(\alpha_1,\dots, \alpha_q)$.
Although expected, this behaviour of the zeros of $L_{\vec{n}}^{\vec{\alpha}}$ was unknown (as far as we know).

\medskip

We can also use well-known properties of the zeros of Laguerre polynomials to prove results on the zeros of generalized Bell polynomials:

\begin{Cor}\label{haf} Assume that $\phi_i\ge 0$. Then
\begin{equation}\label{liz}
-4n-\alpha_n-2<\zeta_1(n)\le\xi_1, \quad n\ge 1,
\end{equation}
where $\xi_1$ denotes the leftmost zero of the Bell polynomial $\Be_n$ and
\begin{equation}\label{aln}
\alpha_n=\max \{-1,\phi_1-1,\phi_2-2,\dots, \phi_n-n\}.
\end{equation}
\end{Cor}

\medskip

When in the sequence $\phi$ there are negative terms, the zeros of the generalized Bell polynomials might lose some of the properties in Theorem \ref{pcgbi} and non real zeros can appear.
Hence, the following problem remains as a challenge.

\medskip

\noindent \textit{Problem 1}: Determine necessary or sufficient conditions for a sequence $\phi=(\phi_i)_i$ of real numbers, with $\phi_i<0$ for some $i\ge 1$, such that
the zeros of the generalized Bell polynomials are all real.
\medskip

The hypergeometric representation (\ref{hypg}) says that Problem 1 is equivalent to determine when the hypergeometric function
$$
\pFq{n}{n}{1+\phi_1,\dots, 1+\phi_n}{\phi_1,\dots, \phi_n}{x}
$$
has only real zeros when $\phi_i<0$ for some $i$.  This case is not include in the Ki and Kim characterization of when $\pFq{n}{n}{a_1,\dots, a_n}{b_1,\dots, b_n}{x}$ has only real zeros \cite{KiKim} (see also \cite{sok}).

We include some partial results on Problem 1. In particular, we extend Theorem \ref{pcgbi} for the case when only one of the $\phi_{i}$'s is negative (Theorem \ref{pcgb2}).

In the last section, we consider the case when $\phi_i \not =0$ only for finitely many $i$'s.
As the main result, we prove that if $\phi_i\not =-1,-2,\dots $, then all the zeros of $\Be_n^\phi$ are reals for $n$ big enough (Theorem \ref{ede}). In that case, the negative zeros of $\Be_{n+1}^\phi$ interlace the negative zeros of $\Be_n^\phi$, and the positive zeros of $\Be_{n+1}^\phi$ interlace the positive zeros of $\Be_n^\phi$.

\section{Basic properties of the generalized Bell polynomials} \label{pre}

We start proving the recurrence (\ref{rrgp}). Actually, we will prove a more general recurrence. In order to do that, we need some notation. For $l\ge 1$, we denote $\phi^{\{ l\}}$ for the sequence
\begin{equation}\label{que}
\phi^{\{ l\} }_i=\begin{cases} \phi_i, &1\le i\le l-1,\\
\phi_{i+1}, &l\le i,\end{cases}
\end{equation}
that is, $\phi^{\{ l\} }$ is que sequence obtained by removing the term $\phi_l$ from $\phi$. We also denote
$$
\Phi_i^{n,l}=\Phi_i^n(\phi^{\{ l\}}).
$$

For a positive integer $l$ and a real number $M$ write $\phi^{l,M}$ for the sequence
\begin{equation}\label{sph}
\phi^{l,M}_i=\phi_i+M\delta_{i,l}
\end{equation}
(where $\delta_{i,l}$ denotes de Kronecker delta).

\begin{Lem}\label{recf}
For all $l\ge 1$ and $n\ge l-1$, we have
\begin{equation}\label{rrgpl}
\Be_{n+1}^\phi(x)=\phi_{l}\Be_n^{\phi^{\{ l\}}}(x)+x\left(1+\frac{d}{dx}\right)\Be_{n}^{\phi^{\{l\}}}(x).
\end{equation}
\end{Lem}

\begin{proof}
From the definition (\ref{dpm}), it is easy to see that, for $n\ge l-1$,
\begin{equation}\label{xxx}
\Phi_i^{n+1}=\begin{cases} \Phi_0^{n+1}, &i=0,\\
\Phi_i^{n,l}+\phi_l\Phi_{i-1}^{n,l}, &1\le i\le n-1,\\
\phi_l\Phi_{n-1}^{n,l},&i= n.\end{cases}
\end{equation}
Hence, an easy computation using (\ref{xxx}) and (\ref{rbel}) gives
\begin{align*}
\Be_{n+1}^\phi(x)&=\sum_{j=0}^{n+1}\Phi_{n+1-j}^{n+1}\Be_j(x)\\
&=\phi_l\Be_{n}^{\phi^{\{l\}}}(x)+\sum_{j=0}^{n}\Phi_{n-j}^{n,l}\Be_{j+1}(x)\\
&=\phi_l\Be_{n}^{\phi^{\{l\}}}(x)+\sum_{j=0}^{n}\Phi_{n-j}^{n,l}x\left(1+\frac{d}{dx}\right)\Be_{j}(x)\\
&=\phi_l\Be_{n}^{\phi^{\{l\}}}(x)+x\left(1+\frac{d}{dx}\right)\left[\sum_{j=0}^{n}\Phi_{n-j}^{n,l}\Be_{j}(x)\right]\\
&=\phi_l\Be_{n}^{\phi^{\{l\}}}(x)+x\left(1+\frac{d}{dx}\right)\Be_{n}^{\phi^{\{l\}}}(x).
\end{align*}

In particular, for $l=n+1$, the identity (\ref{rrgpl}) provides the recurrence relation (\ref{rrgp}) for the generalized Bell polynomials.
\end{proof}

We next prove the identity (\ref{cff}).

Let us remind that the numbers $\rho_{n,j}^\phi$, $0\le j\le n$, are defined from the equation
$$
\prod_{i=1}^n(x+\phi_i)=\sum_{j=0}^n\rho_{n,j}^\phi x(x-1)\cdots (x-j+1).
$$

Using induction on $n$, it is enough to prove that if (\ref{cff}) holds then also $\Be_{n+1}^\phi(z)=\sum_{j=0}^{n+1}\rho_{n+1,j}^\phi x^j$.
From the identity
\begin{align*}
\prod_{i=1}^{n+1}(x+\phi_i)&=\left[\sum_{j=0}^n\rho_{n,j}^\phi x(x-1)\cdots (x-j+1)\right](x+\phi_{n+1})\\
&=\sum_{j=0}^n\rho_{n,j}^\phi x(x-1)\cdots (x-j+1)(x-j+\phi_{n+1}+j),
\end{align*}
we easily get
$$
\rho_{n+1,j}^\phi=\begin{cases} \phi_{n+1}\rho_{n,0}^\phi, & j=0,\\
\rho_{n,j-1}^\phi+(\phi_{n+1}+j)\rho_{n,j}^\phi,& 1\le j\le n,\\
\rho_{n,n}^\phi,&j=n+1. \end{cases}
$$
Hence, using the induction hypothesis and the recurrence relation (\ref{rrgp}), we get
$$
\sum_{j=0}^{n+1}\rho_{n+1,j}^\phi x^j=x\left(1+\frac{d}{dx}\right)\Be_n^\phi(x)+\phi_{n+1}\Be_n^\phi(x)=\Be_{n+1}^\phi(x).
$$

The hypergeometric representation (\ref{hypg}) is now an easy consequence of the identity (\ref{cff}). Indeed, if we write $\rho_{n,l}^\phi=0$, for $l>n$, we have
\begin{align*}
\pFq{n}{n}{1+\phi_1,\dots, 1+\phi_n}{\phi_1,\dots, \phi_n}{x}&=\sum_{j=0}^\infty \frac{\prod_{i=1}^{n}(j+\phi_i)}{j!\, \prod_{i=1}^{n}\phi_i}x^j\\
&=\frac{1}{\prod_{i=1}^{n}\phi_i}\sum_{j=0}^\infty \frac{\sum_{l=0}^n\rho_{n,l}^\phi j(j-1)\cdots (j-l+1)}{j!}x^j\\
&=\frac{1}{\prod_{i=1}^{n}\phi_i}\sum_{j=0}^\infty \frac{\sum_{l=0}^j\rho_{n,l}^\phi j(j-1)\cdots (j-l+1)}{j!}x^j\\
&=\frac{1}{\prod_{i=1}^{n}\phi_i}\sum_{j=0}^\infty\sum_{l=0}^j\frac{\rho_{n,l}^\phi }{(j-l)!}x^j\\
&=\frac{e^x}{\prod_{i=1}^{n}\phi_i}\sum_{l=0}^n\rho_{n,l}^\phi x^l =\frac{e^x}{\prod_{i=1}^{n}\phi_i}\Be_{n}^\phi(x),
\end{align*}
where we assume $\phi_i\not =0,-1,-2,\dots $.

\medskip

As explained in the introduction, the hypergeometric representation (\ref{hypg}) provides a connection between the generalized Bell polynomial and the Laguerre multiple polynomials of the first kind (see \cite[p. 627]{Ism}).
Indeed, the identity (\ref{cgblm}) follows easily by comparing the identity (\ref{hypg}) for the generalized Bell polynomials with the hypergeometric expression for the polynomials $L_{\vec{n}}^{\vec{\alpha}}$ displayed in \cite[p. 628]{Ism}.

\section{Zeros}\label{zer}

In this Section we will prove an stronger version of Theorem \ref{pcgbi}.

\begin{Theo}\label{pcgb} Assume that $\phi_i\ge 0$.
\begin{enumerate}
\item If there exists $i$ such that $\phi_i=0$, writing $i_0=\min \{i:\phi_i=0\}$, we have $\Be_n^\phi(0)\not =0$, $0\le n\le i_0-1$, and $\Be_n^\phi(0) =0$, $i_0\le n$. Otherwise $\Be_n^\phi(0)\not =0$, for all $n$.
\item The polynomial $\Be_n^\phi$ has $n$ simple and non positive zeros, and for all $l\ge 1$, the negative zeros of $\Be_{n+1}^\phi$ interlace the negative zeros of $\Be_{n}^{\phi^{\{l\}}}$.
\item We denote $\zeta_k(n,\phi)$, $1\le k\le n$, for the $k$-th zero of $\Be_n^\phi$, arranging the zeros in increasing order (to simplify the notation and when the context allows it we sometimes will write $\zeta_k$, $\zeta_k(n)$ or $\zeta_k(\phi)$). Then $\zeta_k(\phi)$ is a decreasing function of $\phi$.
\item For a positive integer $l$ and a real number $M\not =0$ with $M>-\phi_l$, consider the sequence $\phi^{l,M}$ (see (\ref{sph})). For $l\le n$, if $M>0$ then the negative zeros of $\Be_{n}^{\phi^{l,M}}$ interlace the negative zeros of $\Be_n^\phi$, and if $M<0$ the negative zeros of $\Be_n^\phi$ interlace the negative zeros of $\Be_{n}^{\phi^{l,M}}$.
\end{enumerate}
\end{Theo}

\begin{proof}

The proof of (1) is trivial, since $\Be_n^\phi(0)=\prod_{i=1}^n\phi_i$.

\medskip

The proof of (2) follows using a similar approach to that for the Bell polynomials (see \cite{har}).

We proceed by induction on $n$.

For $n=1$, we have
$$
\Be_1^\phi(x)=x+\phi_1.
$$
Since $\phi_1\ge 0$, it follows that $\Be_1^\phi$ has a non positive zero.

Assume next that for all sequence of nonnegative numbers $\psi$, $\Be_n^\psi$  has $n$ simple zeros that are non positive real numbers. In particular, the polynomial $\Be_n^{\phi^{\{l\}}}$ has $n$ simple zeros that are non positive real numbers:
$$
\zeta_1<\cdots<\zeta_n\le 0.
$$
Write $r_n(x)=\Be_n^{\phi^{\{l\}}}(x)e^x$, so that the zeros of $r_n$ are the same as those of the polynomial $\Be_n^{\phi^{\{l\}}}$ and with equal multiplicity.
Consider next the function $s_n(x)=xr_n'(x)$. Since $s_n$ has $n$ simple zeros, one of then at $x=0$, $(-1)^ns_n(x)>0$ when $x\to -\infty$, $s_n(x)>0$ when $x\to +\infty$
and the zeros of $r_n$ are simple, then we have
\begin{equation}\label{psc}
s_n(\zeta_i)s_n(\zeta_{i+1})<0,\quad i=1,\dots, n-1, \quad (-1)^n s_n(\zeta_1)>0, s_n(\zeta_n)<0.
\end{equation}
The identities (\ref{rrgpl}) and (\ref{rrgp}), for $n\ge l-1$ and $n\le l-2$, respectively,  give
\begin{equation}\label{mesm}
s_n(x)=\begin{cases}R_{n+1}(x)-\phi_{l}r_n(x),&n\ge l-1,\\ R_{n+1}(x)-\phi_{n+1}r_n(x),&n\le l-2,\end{cases}
\end{equation}
where $R_{n+1}(x)=\Be_{n+1}^{\phi}(x)e^x$ (take into account that for $n\le l-2$ then $\Be_n^{\phi^{\{l\}}}(x)=\Be_n^{\phi}(x)$). Hence, $s_n(\zeta_i)=R_{n+1}(\zeta_i)$ and then (\ref{psc}) gives
$$
R_{n+1}(\zeta_i)R_{n+1}(\zeta_{i+1})<0,\quad i=1,\dots, n-1,\quad (-1)^{n+1} R_{n+1}(\zeta_1)<0, R_{n+1}(\zeta_n)<0.
$$
Since $(-1)^{n+1}R_{n+1}(x)>0$ when $x\to -\infty$ and $R_{n+1}(0)=\phi_1\cdots \phi_{n}\phi_{n+1}\ge 0$, this shows that $R_{n+1}(x)=\Be_{n+1}^{\phi}(x)e^x$ has  $n+1$ real and non positive zeros which interlace those of $r_n(x)=\Be_n^{\phi^{\{l\}}}(x)e^x$. This proves (2).

\medskip

Using the hypergeometric representation (\ref{hypg}), part (2) provides an alternative proof to the fact that for $\phi_i>0$,
the function $\pFq{n}{n}{1+\phi_1,\dots, 1+\phi_n}{\phi_1,\dots, \phi_n}{x}$ has a finite number of zeros and they are all real. For a more general result on zeros of hypergeometric functions see \cite[Theorem 3]{KiKim}.

\medskip

We next prove (3).

Write $\zeta_k=\zeta_k(n,\phi)$ for the $k$-th zero of $\Be_n^\phi$. Since $\Be_n^\phi$ only depends on $\phi_i$, $1\le i\le n$, we have that $\zeta_k$ is an smooth function of each $\phi_i$. In order to prove (3), it is enough to prove that $\partial \zeta_k(\phi)/\partial \phi_i<0$, $1\le i\le n$. To simplify the notation, we write $\zeta=\zeta_k$.
Since $\Be_n^\phi(\zeta)=0$, by deriving with respect to $\phi_i$, we deduce
$$
\frac{d \Be_n^\phi(x)}{dx}\Big|_{x=\zeta}\frac{\partial \zeta(\phi)}{\partial \phi_i}+\frac{\partial \Be_n^\phi(x)}{\partial \phi_i}\Big|_{x=\zeta}=0.
$$
Since $\Be_n^\phi$ has simple and non positive zeros and $\sign(\lim_{x\to -\infty}\Be_n^\phi(x))=(-1)^n$, it follows that
\begin{equation}\label{sig}
\sign \frac{\partial \zeta(\phi)}{\partial \phi_i}=(-1)^{n+k+1}\sign \left(\frac{\partial \Be_n^\phi(x)}{\partial \phi_i}\Big|_{x=\zeta}\right).
\end{equation}
A simple computation shows that (see (\ref{dpm}))
$$
\frac{\partial \Phi^n_l(\phi)}{\partial \phi_i}=\Phi^{n-1}_{l-1}(\phi^{\{i\}}).
$$
Hence, from (\ref{dpbg}), we deduce
\begin{align*}
\frac{\partial \Be_n^\phi(x)}{\partial \phi_i}&=\sum_{j=0}^n\frac{\partial \Phi_{n-j}^n (\phi)}{\partial \phi_i}\Be_j (x)\\
&=\sum_{j=0}^{n-1}\Phi_{n-j}^{n-1} (\phi^{\{i\}})\Be_j(x)=\Be_{n-1}^{\phi^{\{i\}}}(x).
\end{align*}
Using the part (2) of this theorem, we have that
$$
\sign \left(\frac{\partial \Be_n^\phi(x)}{\partial \phi_i}\Big|_{x=\zeta}\right)=\sign \Be_{n-1}^{\phi^{\{i\}}}(\zeta)=(-1)^{n+k}.
$$
Hence, using (\ref{sig}), we finally find
$$
\sign \frac{\partial \zeta(\phi)}{\partial \phi_i}=-1.
$$

\medskip

We finally prove (4). A simple computation gives
$$
\Phi_i^n(\phi^{l,M})=\Phi_i^n(\phi)+M\Phi_{i-1}^{n-1}(\phi^{\{l\}}),
$$
for $1\le i\le n$ and $l\le n$, and so
\begin{equation}\label{mcc}
\Be_{n}^{\phi^{l,M}}(x)=\Be_n^\phi(x)+M\Be_{n-1}^{\phi^{\{l\}}}(x).
\end{equation}
Write $\zeta_i$ for the zeros of $\Be_n^\phi(z)$, so that
$$
\zeta_1<\cdots<\zeta_n\le 0.
$$
Hence $\Be_{n}^{\phi^{l,M}}(\zeta_i)=M\Be_{n-1}^{\phi^{\{l\}}}(\zeta_i)$. It is now enough to take into account the part (2).
\end{proof}

\bigskip

When $\phi_i\ge0$, $i\ge 1$, the zeros of the generalized Bell polynomials $(p_n)_n$  seem to enjoy some other regularities, which already are enjoyed for the zeros of Laguerre polynomials. For instance, in \cite{drjo,drmu}, it is proved that, for $s>0$, the zeros of $L_n^{\alpha+s}$ interlace the zeros of $L_n^\alpha$ for all $n\ge 0$ if and only if $s\in(0,2]$. Taking into account the connection between Laguerre polynomials
and generalized Bell polynomials, we can pose the following problem. For a given sequence of non-negative numbers $\phi=(\phi_i)_i$ and a positive number $s$, define $s+\phi=(s+\phi_i)_i$. Study then the supremum of the set
$$
\mathcal A=\{s>0 : \mbox{the zeros of $\Be_n^{s+\phi}$ interlace the zeros of $\Be_n^\phi$, for all $n\ge 0$}\}.
$$
We have computational evidence which shows that $1\in \mathcal A$, but it would be interesting to calculate $\sup \mathcal A$. This supremum will depend on the sequence $\phi$. Indeed, for $\phi_i=\alpha+i$ (the Laguerre case), we know that $\sup \mathcal A=2$ (and the supremum is actually a maximum), but this is not true in general: for instance, for $\phi=\{1/2,0,0,\dots\}$, we have that the zeros of $\Be_4^{\phi}$ do not interlace the zeros of $\Be_4^{3/2+\phi}$, and so $\sup \mathcal A\le 3/2$.

\bigskip

Take now a multi-index $\vec{n}=(n_1,\dots, n_q)\in \NN^q$ and parameters $\vec{\alpha}=(\alpha_1,\dots, \alpha_q)$,
and denote by $\zeta_k(\vec{\alpha},\vec{n})$, $k=1,\dots,|\vec{n}|=\sum_{i=1}^qn_i$, the zeros of the Laguerre multiple polynomial of the first kind $L_{\vec n}^{\vec \alpha}$ (arranged in increasing order).
According to (\ref{cgblm}), $\zeta_k(\vec{\alpha},\vec{n})=-\zeta_k(|\vec{n}|,\phi^{\vec{\alpha},\vec{n}})$, where the sequence $\phi^{\vec{\alpha},\vec{n}}$ is defined by (\ref{anph}) and $\zeta_k(|\vec{n}|,\phi^{\vec{\alpha},\vec{n}})$ is the $k$-th zero of the generalized Bell polynomial $\Be_{|\vec{n}|}^{\phi^{\vec{\alpha},\vec{n}}}$.
Take other parameters $\vec{\gamma}=(\gamma_1,\dots, \gamma_q)$, such that $\alpha_i\le \gamma_i$. Then $\phi^{\vec{\alpha},{\vec{n}}}\le \phi^{\vec{\gamma},{\vec{n}}}$ (\ref{anph}) and
according to the part (2) of Theorem \ref{pcgb}, $\zeta_k(|\vec{n}|,\phi^{\vec{\gamma},\vec{n}})\le \zeta_k(|\vec{n}|,\phi^{\vec{\alpha},\vec{n}})$. We can then conclude that the zeros of the Laguerre multiple polynomial of the first kind $L_{\vec n}^{\vec \alpha}$ are increasing functions of the parameters $\vec{\alpha}$.

\medskip

We next prove Corollary \ref{haf}.

\begin{proof}[Proof of Corollary \ref{haf}]

For a given sequence of real numbers $\phi_i\ge 0$, $i\ge 1$, on the one hand, the part (3) of Theorem \ref{pcgb} gives that $\zeta_1(n,\phi)\le\xi_1$, where $\xi_1$ denotes the leftmost zero of the Bell polynomial $\Be_n$.
On the other hand, if we take $\alpha_n$ defined by (\ref{aln}) we have
$$
\phi_i=\alpha_n+i+\phi_i-\alpha_n-i\ge \alpha_n+i=\phi ^{\alpha_n}_i.
$$
Since $\Be_n^{\phi^{\alpha_n}}(x)=L_n^{\alpha_n}(-x)$, the part (3) of Theorem \ref{pcgb} gives that $\zeta_1(n,\phi)>-x_{n}$,
where $x_{n}$ denotes the largest zero of the Laguerre polynomial $L_n^{\alpha_n}$. It is then enough to use the well-known estimate $x_{n}<4n+\alpha_n+2$ (see \cite[Theorem 6.31.2]{Sze}).
\end{proof}

\bigskip

In the rest of this section, we include some partial results on the Problem 1.

The following theorem solves completely the Problem 1 when only one of the $\phi_i$'s is negative.

\begin{Theo}\label{pcgb2} Assume that there exists $m\ge 1$ such that $\phi_m<0$ and $\phi_i> 0$, $i\not =m$. Then:
\begin{enumerate}
\item For $1\le n\le m-1$, the polynomial $\Be_n^\phi$ has $n$ simple and negative zeros; for $n\ge m$  the polynomial $\Be_n^\phi$ has $n-1$ simple and negative zeros and one positive zero.
For all $l\ge 1$, $l\not =m$, the negative zeros of $\Be_{n+1}^\phi$ interlace the negative zeros of $\Be_{n}^{\phi^{\{l\}}}$, and for $n\ge m$ the positive zero of $\Be_{n+1}^\phi$ is smaller that the positive zero of $\Be_{n}^{\phi^{\{l\}}}$. For $l=m$, the zeros of $\Be_{n+1}^\phi$ interlace the zeros of $\Be_{n}^{\phi^{\{m\}}}$ (which are all negative).
\item We denote $\zeta_k(n,\phi)$, $1\le k\le n$, for the $k$-th zero of $\Be_n^\phi$, arranging the zeros in increasing order (to simplify the notation and when the context allows it we sometimes will write $\zeta_k$, $\zeta_k(n)$ or $\zeta_k(\phi)$).Then $\zeta_k(n,\phi)$ is a decreasing function of $\phi$, for $k=1,\dots, n-1$, i.e., each negative zero is a decreasing function of $\phi$. The positive zero $\zeta_n(n,\phi)$ is also a decreasing function of the parameter $\phi_m$ but an increasing function of the positive parameters $\phi_i$, $i\not =m$.
\item For a positive integer $l$ and a real number $M$ consider the sequence $\phi^{l,M}$ (see (\ref{sph})). Then, for $n\ge m, l$, we have:
\begin{enumerate}
\item If $l\not =m$ and $M\not =0$ with $M>-\phi_l$, then for $M>0$ the negative zeros of $\Be_{n}^{\phi^{l,M}}$ interlace the negative zeros of $\Be_n^\phi$, and the positive zero of $\Be_{n}^{\phi^{l,M}}$ is greater than the positive zero of $\Be_n^\phi$. For $M<0$ the zeros of $\Be_n^\phi$ interlace the zeros of $\Be_{n}^{\phi^{l,M}}$ and the positive zero of $\Be_n^\phi$ is greater than the positive zero of $\Be_{n}^{\phi^{l,M}}$.
\item If $l=m$ and $M\not =0$, then for $M>0$ the zeros of $\Be_{n}^{\phi^{l,M}}$ interlace the zeros of $\Be_n^\phi$. For $M<0$ the zeros of $\Be_n^\phi$ interlace the zeros of $\Be_{n}^{\phi^{l,M}}$.
\end{enumerate}
\end{enumerate}
\end{Theo}

\begin{proof}

Assume first that $l\not =m$. For $1\le n\le m-1$ the proof is just as that of part (2) in Theorem \ref{pcgb}, but now all the zeros are negative because $\phi_i>0$, $i=1,\dots , m-1$.

Take next $n=m-1$ and denote by
$$
\zeta_1<\cdots<\zeta_{n}< 0
$$
the zeros of $\Be_{n}^{\phi^{\{l\}}}$.
Write $r_n(x)=\Be_n^{\phi^{\{l\}}}(x)e^x$, so that the zeros of $r_n$ are the same of the $\Be_n^{\phi^{\{l\}}}$ and with equal multiplicity.
Consider next the functions $s_n(x)=xr_n'(x)$ and $R_{n+1}(x)=\Be_{n+1}^{\phi}(x)e^x$. The proof follows as in Theorem \ref{pcgb} except for the fact that $R_{n+1}(0)=\phi_1\cdots \phi_{n}\phi_{n+1}<0$, and hence $\Be_n^\phi$ has to have a positive zero.

We next proceed by induction on $n\ge m$. Take, as before,
$$
\zeta_1<\cdots<\zeta_{n-1}<0<\zeta_n
$$
the zeros of the polynomial $\Be_n^{\phi^{\{l\}}}$.
Write $r_n(x)=\Be_n^{\phi^{\{l\}}}(x)e^x$ and $s_n(x)=xr_n'(x)$. Since $s_n$ has simple zeros except at $x=0$ where it might be a double zero,  $(-1)^ns_n(-x)>0$, $s_n(x)>0$ as $x\to +\infty$,
and the zeros of $r_n$ are simple, we have
\begin{equation}\label{psc2}
(-1)^n s_n(\zeta_1)>0,\quad s_n(\zeta_i)s_n(\zeta_{i+1})<0, i=1,\dots, n-2,\quad s_n(\zeta_n)>0,
\end{equation}
and $s_n(\zeta_{n-1})s_n(\zeta_{n})>0$, because $s_n$ has two zeros in $(\zeta_{n-1},\zeta_n)$ (one at $x=0$).
We have (see (\ref{mesm}))
$$
s_n(x)=\begin{cases}R_{n+1}(x)-\phi_{l}r_n(x),&n\ge l-1,\\ R_{n+1}(x)-\phi_{n+1}r_n(x),&n\le l-2,\end{cases}
$$
where $R_{n+1}(x)=\Be_{n+1}^{\phi}(x)e^x$. Hence, $s_n(\zeta_i)=R_{n+1}(\zeta_i)$ and then (\ref{psc2}) gives
$$
(-1)^{n+1} R_{n+1}(\zeta_1)<0,\quad R_{n+1}(\zeta_i)R_{n+1}(\zeta_{i+1})<0,\quad i=1,\dots, n-2,\quad R_{n+1}(\zeta_n)>0.
$$
Since $(-1)^{n+1}R_{n+1}(x)>0$ as $x\to -\infty$, this shows that $R_{n+1}(x)=\Be_{n+1}^{\phi}(x)e^x$ has  $n-1$ simple and negative zeros in $(-\infty,\zeta_{n-1})$ which interlace the first $n-1$ negative zeros of $r_n(x)=\Be_n^{\phi^{\{l\}}}(x)e^x$.
We also have $R_{n+1}(\zeta_{n-1})=s_n(\zeta_{n-1})=(-1)^{n+n-1+1}=1$, $R_{n+1}(\zeta_{n})=s_n(\zeta_{n})=1$ and $R_{n+1}(0)=\phi_1\cdots \phi_{n}\phi_{n+1}<0$. This implies that $R_{n+1}$ has other negative zero in $(\zeta_n,0)$
and a positive one in $(0,\zeta _n)$. This proves (1) for $l\not =m$.

The case $l=m$ can be proved similarly.

In the proof of part (2), we proceed similarly as in part (3) of Theorem \ref{pcgb}, but now we have to use
$$
\sign \left(\frac{d\Be_n^\phi(x)}{dx}\Big|_{x=\zeta_k}\right)=(-1)^{n+k},\quad 1\le k\le n
$$
(as in Theorem \ref{pcgb}), and
$$
\sign \left(\frac{\partial \Be_n^\phi(x)}{\partial \phi_i}\Big|_{x=\zeta_k}\right)=\sign \Be_{n-1}^{\phi^{\{i\}}}(\zeta_k)=\begin{cases}(-1)^{n+k}, & 1\le i\le n, 1\le k\le n-1,\\
1,& k=n, i=m,\\ -1,& k=n, 1\le i\le n, i\not =m \end{cases}
$$
(because of part (1) above in this theorem).

The proof of part (3) is as in Theorem \ref{pcgb}.
\end{proof}

When two of the $\phi_i$'s are negative then Theorem \ref{pcgb} is no longer true because complex zeros can appear. That is the case, for instance, when
$$
\phi_i=\begin{cases} -m,&i=1,2,\\0,&i\ge 3,\end{cases}
$$
and $m$ a positive integer bigger than $1$: we prove in Lemma \ref{um1} below that $\Be_n^\phi$, $n\ge 2$, has always two non real zeros and $n-3$ negative zeros plus a zero at $x=0$.

\bigskip

Let us remark that using the hypergeometric representation (\ref{hypg}), we can rewrite Theorem \ref{pcgb2} in terms
of the zeros of the hypergeometric function
$$
\pFq{n}{n}{1+\phi_1,\dots,1+\phi_n}{\phi_1,\dots,\phi_n}{z}.
$$
In particular when there exists $m\ge 1$ such that $\phi_m<0$ and $\phi_i> 0$, $i\not =m$, that hypergeometric function has exactly $n$ real zeros, and for $n\ge m$, $n-1$ of them are negative
(and they are decreasing functions of the parameters $\phi_i$) and one is positive (which it is also a decreasing function of the parameter $\phi_m$ but an increasing function of the positive parameters $\phi_i$, $i\not =m$).

\section{The case when $\phi_i \not =0$ only for finitely many $i$'s.}\label{ult}

We next study the case when only finitely many of the $\phi_i$'s are different to zero. Hence assume that
\begin{equation}\label{acd}
\mbox{there exists a positive integer $K$ such that
$\phi_i=0$, for $i\ge K+1$.}
\end{equation}
Write $P$ for the polynomial
\begin{equation}\label{pop}
P(x)=\prod_{i=1}^K(x+\phi_i)=\sum_{j=0}^K\gamma_{j}x^{K-j}.
\end{equation}
Hence, for $n\ge K$,
$$
\Be_n^\phi(x)=\sum_{j=0}^K \gamma_j\Be_{n-j}(x).
$$
Since $\phi_{K+1}=0$, we know that $\Be_n^\phi(0)=0$, $n\ge K+1$.
We next characterize when $x=0$ has multiplicity bigger than one.

\begin{Lem}\label{xzero} Assume (\ref{acd}). The following are equivalent.
\begin{enumerate}
\item There exists a positive integer $l>1$ such that for all $n\ge K+1$, $\Be_{n}^\phi$ has a zero of multiplicity $l$ at $x=0$.
\item There exists a positive integer $l>1$ such that for some $n\ge K+1$, $\Be_{n}^\phi$ has a zero of multiplicity $l$ at $x=0$.
\item There exist $i_1,\dots,i_{l-1}$ such that $\phi_{i_j}=-j$, $j=1,\dots , l-1$.
\end{enumerate}
\end{Lem}

\begin{proof}
Since $\Be_n^{(m)}(0)=m!\, S(n,m)$, using (\ref{stir}), (\ref{dpbg}) and (\ref{pop}), we have, for $n\ge K+1$,
\begin{align*}
(\Be_n^\phi)^{(m)}(0)&=\sum_{j=0}^K\gamma_j\Be_{n-j}^{(m)}(0)=m!\sum_{j=0}^K\gamma_jS(n-j,m)=m!\, \sum_{j=0}^K\gamma_j\sum_{i=0}^m\frac{(-1)^{m-i}i^{n-j}}{(m-i)!\, i!}\\
&=m!\sum_{i=0}^m \frac{(-1)^{m-i}i^{n-K}}{(m-i)!\, i!}\sum_{j=0}^K\gamma_ji^{K-j}\\
&=m!\sum_{i=0}^m \frac{(-1)^{m-i}i^{n-K}}{(m-i)!\, i!}P(i)
\\
&=m!\sum_{i=0}^m \frac{(-1)^{m-i}i^{n-K}}{(m-i)!\, i!}\prod_{j=1}^K(i+\phi_j).
\end{align*}
The proof now follows straightforwardly.
\end{proof}

From (\ref{rbel}) it follows that, for $n\ge K$,
\begin{equation}\label{idc}
\Be_{n+1}^\phi(x)=x\left(1+\frac{d}{dx}\right)\Be_n^\phi(x).
\end{equation}
This fact has important consequences on the zeros of the polynomials $\Be_n^\phi(x)$, $n\ge K$. Indeed, using \cite[Lemma 8]{dur}, we have the following.

\begin{Cor}\label{mie2}
If we assume (\ref{acd}), then:
\begin{enumerate}
\item The number of real zeros of the polynomial $\Be_n^\phi$, $n\ge K$, is at least $n-K+1$ if $K$ is odd, and $n-K$ if $K$ is even.
\item If there exists $n_0\ge K$ such that $\Be_{n_0}^\phi$ has only real zeros, then $\Be_n^\phi$ has only real zeros for $n\ge n_0$.
\item If there exists $n_0\ge K$ such that $\Be_{n_0}^\phi$ has only real and simple zeros, then $\Be_n^\phi$ has only real and simple zeros for $n\ge n_0$.
\end{enumerate}
\end{Cor}

We already know that if $\phi_i<0$ for some $i$, positive zeros can appear. It was proved in \cite[Section 4]{dur} that assuming
\begin{equation}\label{asu}
\phi_i\not =-1,-2,\dots,\quad 1\le i\le K,
\end{equation}
for $n$ big enough, the exact number of positive zeros of $p_n$ does not depend on $n$ and it is given by the number of elements $s$ of the set
\begin{equation}\label{eha}
\Hh=\{l: \mbox{$l$ is a positive integer and $P(l)P(l+1)<0$}\}.
\end{equation}

In particular, under the assumptions (\ref{acd}) and (\ref{asu}), according to Lemma \ref{xzero} the polynomial $\Be_n^\phi$, $n\ge K+1$, has a simple zero at $x=0$.

Using the identity (\ref{idc}), we have the following.

\begin{Lem}\label{illa} Assume (\ref{acd}) and (\ref{asu}) and let $s$ be the number of elements of $\Hh$. Then there exists $n_0$ such that, for $n\ge n_0$, $\Be_n^\phi$ has exactly $s$ positive (and simple) zeros, and if for some $n_1\ge n_0$, $\Be_{n_1}^\phi$ has real and simple zeros, then for all
$n\ge n_1$, the zeros of $\Be_n^\phi$ are also real and simple. Moreover, in that case, the $n-s$ negative zeros of $\Be_{n+1}^\phi$ interlace the $n-s-1$ negative zeros of $\Be_{n}^\phi$, and the $s$ positive zeros of $\Be_{n+1}^\phi$ interlace the $s$ positive zeros of $\Be_{n}^\phi$.
\end{Lem}

\begin{proof}

The first part is a consequence of \cite[Theorem 9]{dur} and Corollary \ref{mie2}.

For the second part, write $\zeta_{j}^-$, $1\le j\le n-s-1$, for the negative zeros of $\Be_n^\phi $ arranged in decreasing order, and $\zeta_j^+$, $1\le j\le s$, for the positive zeros arranged in increasing order.

The identity (\ref{idc}) implies that
\begin{equation}\label{hos}
\Be_{n+1}^\phi(\xi_j^\pm)=\xi_j^\pm(\Be_{n}^\phi)'(\xi_j^\pm).
\end{equation}
Hence, in each interval $(\zeta_{j+1}^-,\zeta_j^-)$ and $(\zeta_{j}^+,\zeta_{j+1}^+)$, $\Be_{n+1}^\phi$ has an odd number of zeros. Since $\Be_{n+1}^\phi$ has $s$ positive zeros,
(\ref{hos}) also implies that $\Be_{n+1}^\phi$ can not have any zero in $(\zeta_s^+,+\infty)$, and hence we conclude that $\Be_{n+1}^\phi$ has a zero in $(0,\zeta_1^+)$. (\ref{hos}) then shows that $\Be_{n+1}^\phi$ has also a zero in $(\zeta_1^-,0)$. This proves the lemma.
\end{proof}

We next prove that assuming  (\ref{acd}) and (\ref{asu}), for $n$ big enough $\Be_n^\phi$ has always simple and real zeros.

\begin{Theo}\label{ede} If we assume (\ref{acd}) and (\ref{asu}), then there exists $n_0$ such that for $n\ge n_0$, the zeros of $\Be_n^\phi $ are real and simple. Moreover the negative zeros of $\Be_{n+1}^\phi $ interlace the negative zeros of $\Be_n^\phi $, and the positive zeros of $\Be_{n+1}^\phi $ interlace the positive zeros of $\Be_n^\phi $.
\end{Theo}

\begin{proof}

We proceed by induction on $K$.

For $K=1$, the result is a direct consequence of the Obreshkov theorem applied to $\Be_n/x$, $\Be_{n-1}/x$.

\begin{TheoA}\label{obre}
If $p$ and $q$ are real polynomials of degrees differing by $1$ at most and having no common zero, then a necessary and sufficient condition for their zeros to be real and distinct and to interlace each other is that $\lambda p(z)+\mu q(z)$ have only distinct real zeros for all real $\lambda$ and $\mu$.
\end{TheoA}

(see \cite[p. 9]{hou}).

%\begin{TheoA}\label{obre}
%Let $p$ and $q$ be polynomials with only real and simple zeros and with $\deg (p)=\deg(q)+1$. Then the following are equivalent.
%\begin{enumerate}
%\item The zeros of $p$ interlace the zeros of $q$.
%\item For each real numbers $\alpha,\beta$, the polynomial $\alpha p+\beta q$ has only real and simple zeros.
%\end{enumerate}
%\end{TheoA}

Assume next that the theorem is true for $K$.

Let $\phi_i$ be, $1\le i\le K+1$, real numbers with $\phi_i\not =-1,-2,\dots$, for $1\le i\le K+1$. In order to simplify the notation, we write $p_n=\Be_n^{\phi}$ and $q_n=\Be_n^{\phi^{\{K+1\}}}$ (see (\ref{que})).
Since $\Be_n^{\phi^{K+1,-\phi_{K+1}}}=\Be_n^{\phi^{\{K+1\}}}=q_n$ (see (\ref{sph})), the following identity follows from (\ref{mcc}):
\begin{equation}\label{pmi}
q_n(z)=p_n(x)-\phi_{K+1}q_{n-1}(x).
\end{equation}
The induction hypothesis says that the zeros of $q_n$ are simple and real for $n$ big enough. Write $\zeta_j^+$, $1\le j\le s$, for the positive zeros of $q_n$ arranged in increasing order, and $\zeta_j^-$, $1\le j\le n-s-1$, for the negative zeros arranged in decreasing order. Setting $P(x)=\prod_{i=1}^K(x+\phi_i)$, we have that $s$ is the number of elements of the set $\Hh$ (\ref{eha}).
Hence
$$
\zeta^-_{n-s-1}<\dots<\zeta^-_1<0<\zeta^+_1<\dots <\zeta^+_s.
$$
The identity (\ref{pmi}) then gives
\begin{equation}\label{esm}
p_n(\zeta_i^\pm)=\phi_{K+1}q_{n-1}(\zeta_i^\pm).
\end{equation}
Lemma \ref{illa} says that
\begin{align}\label{esm1}
\sign[q_{n-1}(\zeta_{n-i-1}^-)]&=(-1)^{n-1+s-i},\quad s\le i\le n-2,\\\label{esm2}
\sign[q_{n-1}(\zeta_{i}^+)]&=(-1)^{s-i+1},\quad 1\le i\le s.
\end{align}
Hence, on the one hand, (\ref{esm}) and (\ref{esm1}) say that
\begin{align}\label{ton1}
&\mbox{$p_n$ has an odd number of zeros in}\\ \nonumber
&\hspace{3cm}\mbox{each interval $(\zeta_{i+1}^-,\zeta_i^-)$, $1\le i\le n-s-2$},
\end{align}
and so $p_n$ has at least $n-s-2$ negative zeros. On the other hand,
(\ref{esm}) and (\ref{esm1}) say that
\begin{equation}\label{ton2}
\mbox{$p_n$ has an odd number of zeros in each interval $(\zeta_{i}^+,\zeta_{i+1}^+)$, $1\le i\le s-1$},
\end{equation}
that is, $p_n$ has at least $s-1$ positive zeros.
Since $p_n(0)=0$, we have that $p_n$ has at least $n-2$ real zeros.

If $\phi_{K+1}>0$, on the one hand, using (\ref{esm1}) we deduce that $p_n$ has one more zero in the interval $(-\infty,\zeta_{n-s-1}^-)$. On the other hand, using (\ref{esm2}) we deduce that $p_n$ has also a zero in the interval $(\zeta_{s}^+,+\infty)$. Hence $p_n$ has $n$ real zeros, and since they are located in disjoint interval (see (\ref{ton1}) and (\ref{ton2})), we conclude that the zeros of $p_n$ have to be simple.

If $\phi_{K+1}<0$, write $m_0=\lfloor -\phi_{K+1} \rfloor$, so that $m_0<-\phi_{K+1}<m_0+1$ ($\lfloor x\rfloor$ denotes the floor function).
Setting $\tilde P(x)=\prod_{i=1}^{K+1}(x+\phi_i)$,
$$
\tilde \Hh=\{l: \mbox{$l$ is a positive integer and $\tilde P(l)\tilde P(l+1)<0$}\},
$$
and denoting by $\tilde s$ the number of elements of $\tilde \Hh$, we have that $\tilde P(x)=(x+\phi_{K+1})P(x)$, and for $n$ big enough $\tilde s$ is the number of positive zeros of $p_n$.

If $P(m_0)P(m_0+1)>0$, we have that $\tilde P(m_0)\tilde P(m_0+1)<0$ and then $\tilde s=s+1$. That is, $p_n$ has to have $s+1$ positive zeros, and we can then conclude that $p_n$ has to have $n$ real zeros. According to \cite[Theorem 9]{dur}, for $n$ big enough, the positive zeros has to be simple.
Moreover, we have also proved that $p_n$ has $n-s-2$ negative zeros, and then (\ref{ton1}) says that they have to be simple. Therefore, we conclude that all the zeros of $p_n$ have to be real and simple.

Assume next that $P(m_0)P(m_0+1)<0$. We then have $\tilde P(m_0)\tilde P(m_0+1)>0$, and so
$$
\tilde \Hh=\Hh\setminus\{m_0\}.
$$
As a consequence, we deduce that $\tilde s=s-1$. Write
$$
\G=\{1,2,3,\dots\}\setminus \Hh=\{g_i:1\le i\},\quad \tilde \G=\{1,2,3,\dots\}\setminus \tilde \Hh=\{\tilde g_i:1\le i\},
$$
with $g_i<g_j$, $\tilde g_i<\tilde g_j$ if $i<j$. It is then easy to see that $\tilde g_{m_0}=g_{m_0}-1$.

Write $\tilde \zeta_j^-$ for the negative zeros of $p_n$ arranged in decreasing order, where either
$1\le j\le n-s$ or $1\le j\le n-s-2$. \cite[Theorem 9]{dur} says that for any positive integer $m$ and $n$ big enough, $\tilde \zeta_j^-$, $j=1,\dots ,m$, are simple and they are located in the interval $(v_{m,n},0)$, where
$$
v_{m,n}=-(1+\epsilon)\tilde g_m\frac{\tilde P(\tilde g_m)}{\tilde P(\tilde g_{m}+1)}\left(\frac{\tilde g_m}{\tilde g_{m}+1}\right)^{n-K-2}.
$$
If we take $m=m_0$,  \cite[Theorem 9]{dur} applied to $q_n$ gives that $\zeta_{m_0}^-<u_{n}$, where
$$
u_{n}=-(1-\epsilon) g_{m_0}\frac{ P( g_{m_0})}{ P( g_{m_0}+1)}\left(\frac{ g_{m_0}}{g_{m_0}+1}\right)^{n-K-1}.
$$
We next assume the following:

\medskip
\noindent
\textit{Claim}: $u_n<v_{m_0,n}$.

Hence, on the one hand, $p_n$ has $s-1$  simple and positive zeros, plus a simple zero at $x=0$. This gives that $p_n$ has to have at most $n-s$ negative zeros.
The claim implies that $p_n$ has at least $m_0$ simple and negative zeros in $(v_{m_0,n},0)$, and, since the number of negative zeros is at most $n-s$, we deduce from (\ref{ton1}) that $p_n$ has exactly $n-m_0-s-1$ simple zeros in $(\zeta_{n-s-1}^-,\zeta_{m_0}^-)$. Hence $p_n$ has $n-s$ negative zeros and there is only one negative zero to be located. Since $\phi_{K+1}<0$, using (\ref{esm1}) we deduce that $p_n$ has an even number of zeros in the interval $(-\infty,\zeta_{n-s-1}^-)$, hence $p_n$ does not vanish in $(-\infty,\zeta_{n-s-1}^-)$. This implies that $p_n$ has exactly $m_0+1$ simple and negative zeros in $(v_{m_0,n},0)$. That is, $p_n$ has $n$ simple and real zeros.

We finish proving the claim. Since $\tilde g_{m_0}=g_{m_0}-1$, we have to prove that
\begin{align*}
(1-\epsilon) g_{m_0}\frac{ P( g_{m_0})}{ P( g_{m_0}+1)}&\left(\frac{ g_{m_0}}{g_{m_0}+1}\right)^{n-K-1}\\ &>
(1+\epsilon)( g_{m_0}-1)\frac{\tilde P( g_{m_0}-1)}{\tilde P(g_{m_0})}\left(\frac{ g_{m_0}-1}{g_{m_0}}\right)^{n-K-2}.
\end{align*}
On the one hand, since $g_{m_0}\in \G$, we get $P( g_{m_0})/P( g_{m_0}+1)>0$. And, on the other hand, since $g_{m_0}-1=\tilde g_{m_0}\in \tilde \G$ we have $\tilde P( g_{m_0}-1)/\tilde P(g_{m_0})>0$. Hence, the previous inequality follows for $n$ big enough (depending on $m_0$) because
$$
\frac{ g_{m_0}}{g_{m_0}+1}>\frac{ g_{m_0}-1}{g_{m_0}}.
$$

\end{proof}

\bigskip If we remove the assumption (\ref{asu}) then, in general, we can not guarantee that $\Be_n^\phi$ has only real zeros for $n$ big enough, as the following example shows.

\begin{Lem}\label{um1}
If $\phi_1=\phi_2= -m$, $\phi_i=0$, $i\ge 3$, with $m\in \NN$ and $m\ge 2$, then for all $n\ge 3$ the generalized Bell polynomial $\Be_n^\phi(x)$ vanish at $x=0$, has two non real zeros and $n-3$
negative zeros. Moreover, the negative zeros of $\Be_n^\phi$, $n\ge 3$, interlace the negative zeros of $\Be_{n-1}^\phi$.
\end{Lem}

\begin{proof}

We already know from Corollary \ref{mie2} that the generalized Bell polynomial $\Be_n^\phi(x)$ vanish at $x=0$, and it has at least $n-3$ real zeros.

In this case, for $n\ge 2$, we have
$$
\Be_n^\phi(x)=\Be_n(x)-2m\Be_{n-1}(x)+m^2\Be_{n-2}(x).
$$
Hence, if we write $a_{n,j}$ for the coefficient of $x^j$ of $\Be_n^\phi$, using that
\begin{equation}\label{rsn}
S(n,j)=jS(n-1,j)+S(n-1,j-1)
\end{equation}
(see \cite[p. 266]{GKP}), we get
\begin{align*}
a_{n,j}&=S(n,j)-2m S(n-1,j)+m^2S(n-2,j)
\\&=(j-2m)S(n-1,j)+S(n-1,j-1)+m^2S(n-2,j).
\end{align*}
This shows that $a_{n,j}>0$ for $j\ge 2m $.

Using the asymptotic
\begin{equation}\label{aanx}
S(n,j)\sim \frac{j^n}{j!}
\end{equation}
for fixed $j$ and $n$ big enough (see \cite[p. 121]{bon}), we get that for $1\le j<2m$, $j\not =m$,
\begin{equation}\label{aan}
a_{n,j}\sim \frac{j^{n-2}}{j!}(j-m)^2,
\end{equation}
and so, we can conclude that for $j\not =m$ and $n$ big enough then $a_{n,j}>0$.

For $j=m$, using (\ref{rsn}) and (\ref{aanx}), we have
\begin{equation}\label{aany}
a_{n,m}=S(n-1,m-1)-mS(n-2,m-1)\sim -\frac{(m-1)^{n-2}}{(m-1)!}.
\end{equation}
Hence, we deduce that $a_{n,m}<0$ for $n$ big enough.

Descartes' rule of signs says then that $\Be_n^\phi(x)/x$ has at most $n-3$ negative zeros and at most $2$ positive zeros.

Since $a_{n,j}>0$, $j\not =m$, for $x\ge 0$, we have
\begin{equation}\label{aan2}
\Be_n^\phi(x)\ge x^{m-1}(a_{n,m+1}x^2+a_{n,m}x+a_{n,m-1}).
\end{equation}
Using (\ref{aan}) and (\ref{aany}), we deduce that
$$
a_{n,m}\sim -\frac{(m-1)^{n-2}}{(m-1)!},\quad a_{n,m\pm 1}\sim \frac{(m\pm 1)^{n-2}}{(m\pm 1)!}.
$$
Hence for $n$ big enough
$$
a_{n,m}^2-4a_{n,m+1}a_{n,m-1}\sim \frac{(m-1)^{2n-4}}{(m-1)!^2}\left[1-\frac{4}{m(m+1)}\left(\frac{m+1}{m-1}\right)^{n-2}\right]<0.
$$
And so, from (\ref{aan2}) we deduce that for $n$ big enough $\Be_n^\phi(x)>0$, $x\ge 0$. Then for $n$ big enough, $\Be_n^\phi(x)$ can not have positive zeros and so $\Be_n^\phi(x)$ has exactly $n-2$ real zeros, one at $x=0$ and the other zeros are negative.

Since, for $n\ge 2$,
$$
\Be_n^\phi(x)=x\Be_{n-1}^\phi (x)+\frac{d}{dx}\Be_{n-1}^\phi(x),
$$
it is easy to prove that between two simple and consecutive zeros of $\Be_{n-1}^\phi$ of the same sign $\Be_n^\phi$ has to vanish. This shows that for all $n\ge 2$,  $\Be_n^\phi(x)$ has exactly $n-2$ real zeros, one at $x=0$ and the other zeros are negative, and the negative zeros of $\Be_n^\phi$, $n\ge 3$, interlace the negative zeros of $\Be_{n-1}^\phi$.
\end{proof}

\bigskip

Given real numbers $\gamma_j$, $0\le j\le K$, with $\gamma_0=1$, we write
\begin{equation}\label{ppnxx}
p_n(x)=\sum_{j=0}^K \gamma_j\Be_{n-j},\quad n\ge K.
\end{equation}
Define the polynomial
\begin{equation}\label{popxx}
P(x)=\sum_{j=0}^K\gamma_{j}x^{K-j}.
\end{equation}
If $P$ has real zeros $\theta_i$, $i=1,\dots, K$, then it is clear that $p_n=\Be_n^\phi$, where $\phi_i=-\theta_i$, $1\le i\le K$, and $\phi_i=0$, $i\ge K+1$. Hence, Theorem \ref{ede} says that if
$\theta_i\not=1,2,\dots$, then $p_n$ has only real zeros for $n$ big enough.
We guess that if $P$ satisfies $\theta_i\not=1,2,\dots$, even if $P$ has non-real zeros, the polynomial $p_n$ will have only real zeros for $n$ big enough.
\medskip

\noindent
\textit{Conjecture.} If $P(m)\not =0$ for $m=1,2,\dots $, then for $n$ big enough all the zeros of $p_n$ (\ref{ppnxx}) are real.

\bigskip

It is easy to check that the operator $T$ acting in the linear space of polynomials and defined by $T(p)(x)=x\left(p(x)+p'(x)\right)$ is a real zero increasing operator, in the sense that
for all polynomial $p$ the number of real zeros of $T(p)$ is greater than the number of real zeros of $p$ (see \cite[Lemma 8]{dur}). The conjecture shows that this operator $T$ is likely a real zero increasing operator in a more deep way. Indeed, take a monic polynomial $p$ of degree $K$, and write
$$
p(x)=\sum_{j=0}^K\gamma_j\Be_{n-j}(x)
$$
(so that $\gamma_0=1$). It is easy to check that $p_n=T^n(p)$ (\ref{ppnxx}).

If the polynomial $P(x)=\sum_{j=0}^K\gamma_{j}x^{K-j}$ (\ref{popxx}) satisfies $P(m)\not =0$ for $m=1,2,\dots $, the conjecture says that for $n$ big enough the polynomial $T^n(p)$ has only real zeros, even if all the zeros of $p$ are non-real.

\medskip

%%%%%%%%%%%%%%%%%%%%%%%%%%%%%%%%%%%%%%%%%%%%%%%%%%%%%%%%%%%%%%%%%%%%%%%

%%%%%%%%%%%%%%%%%%%%%%%%%%%%%%%%%%%%%%%%%%%%%%%%%%%%%%%%%%%%%%%%%%%%%%%

\end{document}